\documentclass[11pt]{amsart}
\usepackage[headings]{fullpage}
\usepackage{amssymb}
\usepackage{graphicx} 
\usepackage{enumerate}
\usepackage{mathtools}
\usepackage{color,soul}
\usepackage{cite}
\usepackage{tikz}
\usetikzlibrary{matrix}
\usepackage{caption}
\usepackage{amsmath, amsthm}
\usepackage{comment}
 \usepackage{relsize}
 \usetikzlibrary{cd}
 \usetikzlibrary{decorations.pathreplacing}
 \usepackage{tikz,calc}
\usepackage{color}
\usepackage[margin=1.2in]{geometry}
\usepackage{url}

\usepackage{multicol}

\newtheorem{thm}{Theorem}[section]

\newtheorem{cor}[thm]{Corollary}

\theoremstyle{definition}
\newtheorem{defn}[thm]{Definition}
\newtheorem{ex}[thm]{Example}

\newcommand{\theoremname}{Theorem:}

\usepackage{bbm}

\newcommand{\R}{\mathbb{R}}

\newcommand{\N}{Niebrzydowski }

\begin{document}

\title{{An Invariant of Virtual Trivalent Spatial Graphs}}

\author{Evan Carr}
\address{Elon University}
\email{ecarr10@elon.edu}
\urladdr{}

\author{Nancy Scherich}
\address{Elon University}
\email{nscherich@elon.edu}
\urladdr{http://www.nancyscherich.com}

\author{Sherilyn Tamagawa}
\address{Davidson College}
\email{shtamagawa@davidson.edu}
\urladdr{http://sherilyntamagawa.com}


\thanks{}

\begin{abstract}
 We create an invariant of virtual $Y$-oriented trivalent spatial graphs using colorings by \emph{virtual \N algebras}. This paper generalizes the color invariants using \emph{virtual tribrackets} and \emph{\N algebras} by Nelson-Pico, and  Graves-Nelson-T. We computed all tribrackets, \N algebras and virtual \N algebras of orders 3 and 4, and provide generative code for all data sets.
 
\end{abstract}
\maketitle

\section{Introduction}


A \emph{Y-oriented trivalent spatial graph} is an oriented finite graph embedded in $\R^3$ and which contains only degree 3 vertices.
The $Y$-orientation forbids sources or sinks at the vertices, but requires vertex orientations with \emph{two-in-one-out} or \emph{two-out-one-in}.
A \emph{virtual} $Y$-oriented trivalent spatial graph is a $Y$-oriented trivalent spatial graph with an additional type of crossing, a virtual crossing. Virtual crossings were originally defined for virtual knots in \cite{K}, and we denote virtual crossings by a circle over the crossing in the diagram.
Two virtual $Y$-oriented trivalent spatial graphs are pictured below.

\begin{figure}[h]
\includegraphics[width=240\unitlength]{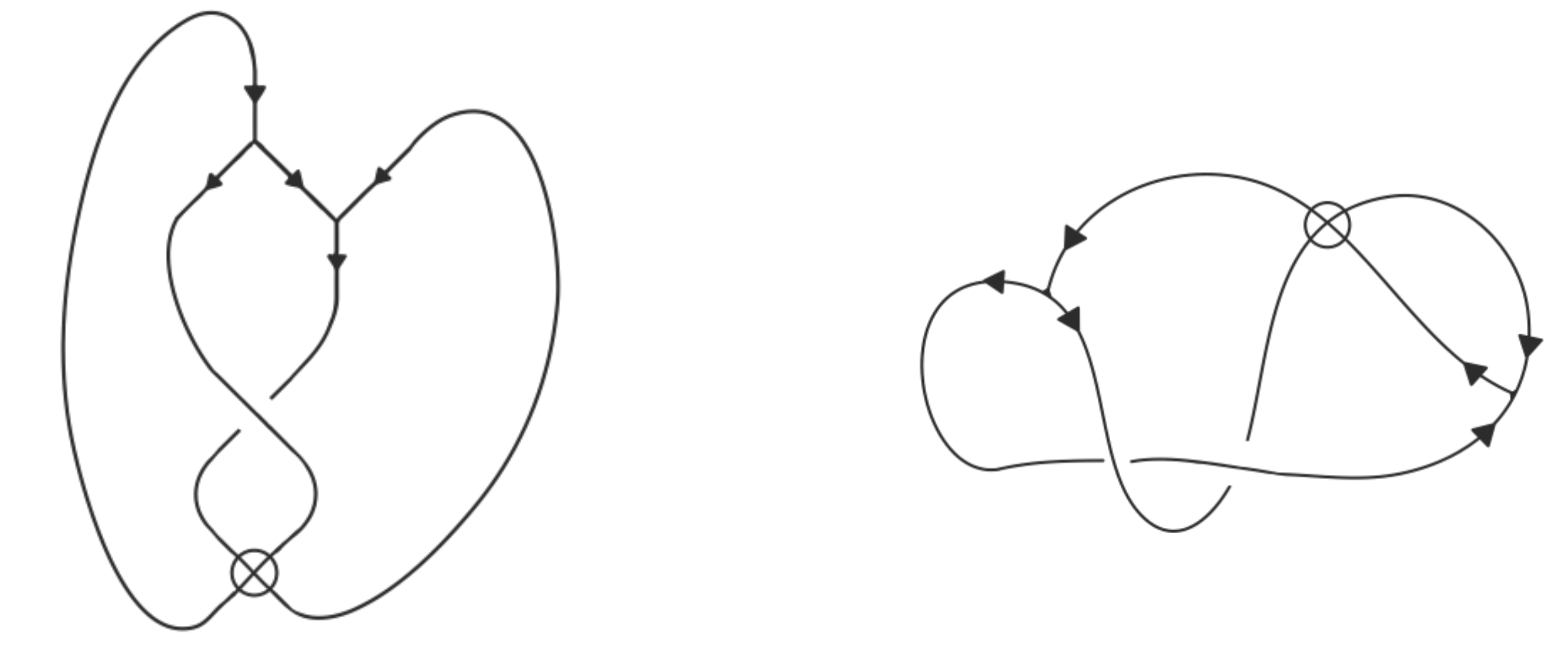}
\caption{Two virtual Y-oriented trivalent spatial graphs.}
\end{figure}

Virtual trivalent graphs have become increasingly important in finite type invariant theory since Bar-Natan and Dancso developed \emph{w-Foams} to study certain types of knotted surfaces in $\R^4$ \cite{BD2}. Murakami, Ohtsuki, and Yamada proved that the Reshetikhin-Turaev invariant of a link in $\R^3$ can be computed using (classical) oriented trivalent graphs \cite{MOY}. 
Trivalent vertices are particularly useful in algebraic computations by labeling the edges of a trivalent graph. For example,  the trivalent vertex orientated with \emph{two-in-one-out} can be used as a visual representation of a Lie bracket, as was done in \cite{BDS}.
The two inputs label the inward oriented arcs of the vertex, and the output of the Lie bracket labels the outward oriented arc. 
Similarly, the vertex oriented with \emph{two-out-one-in} can be used to represent a coproduct.
In this paper, instead of labeling the arcs of a graph, we label/color the planar complement of the graph to attain an invariant.

Coloring invariants of knotted objects have been widely studied in various different contexts. 
For example, \N used an algebraic structure called a \emph{ternary quasigroup} to color planar knot complements \cite{N}. 
Nelson and Pico used \emph{virtual tribrackets} to color planar complements of oriented virtual link diagrams \cite{NP}. 
Extending the ideas of \N and Nelson-Pico,  Graves-Nelson-T. introduced \emph{\N algebras} to color planar complements of $Y$-oriented trivalent spatial graphs \cite{GNT}.

In the original definition of an \N algebra given in \cite{GNT}, there is a well-definedness issue related to the partial product structure. 
Partial products can lead to violations of the Reidemeister moves, which nullifies the use of an \N algebra to define a graph invariant. 
We offer a correction to this issue with the addition of an extra axiom to the definition of an \N algebra, as described in Section \ref{sec:welld}.

In this paper, we combine the ideas  of a virtual tribracket and a \N algebra to define a \emph{virtual \N algebra}, as stated in Definition \ref{def:vnalg} in Section \ref{sec:nalgs}. 
Virtual \N algebras can be used to color the planar complements of virtual $Y$-oriented trivalent spatial graphs. Our main result proves that the number of colorings of a graph by a virtual \N algebra is an invariant of the graph.\\

\noindent \textbf{Main result}(Theorem 3.6) \emph{Let $X$ be a virtual \N algebra and $\Gamma$ a virtual $Y$-oriented trivalent spatial graph. The number of $X$-colorings of the planar complement of $\Gamma$ is invariant under the virtual $Y$-oriented Reidemeister moves, and is therefore an integer-valued invariant of virtual $Y$-oriented trivalent spatial graphs.}\\

\noindent\emph{Organization of the paper.} In Section \ref{sec:Rmoves}, we discuss the virtual $Y$-oriented Reidemeister moves for virtual $Y$-oriented trivalent spatial graphs. Section \ref{sec:nalgs} reviews the foundational definitions of tribrackets, virtual tribrackets, \N algebras, and the associated coloring invariants for these structures. We state our new definition of virtual \N algebras and prove our main result. In Section \ref{sec:welld}, we discuss the well-definedness issues of the original definition of an \N algebra and give a corrected definition of a \emph{partially defined \N algebra}, and extend this notion to \emph{partially defined virtual \N algebras} as well.  Much of the work in this paper is computational, and in Section \ref{sec:comp} we share the computational results and give links to generative code to use the data sets.\\

\noindent \emph{Acknowledgments.} The authors would like to thank Sam Nelson for helpful conversations, and our anonymous referee for catching an important error. The second author was partially supported by the NSERC grant RGPIN-2018-04350. This project occurred in part as an undergraduate student research project through Elon University where the first author was mentored by the second author.

\section{Virtual Y-oriented Reidemeister Moves}\label{sec:Rmoves}

\begin{figure}
    \centering
    \includegraphics[width=420\unitlength]{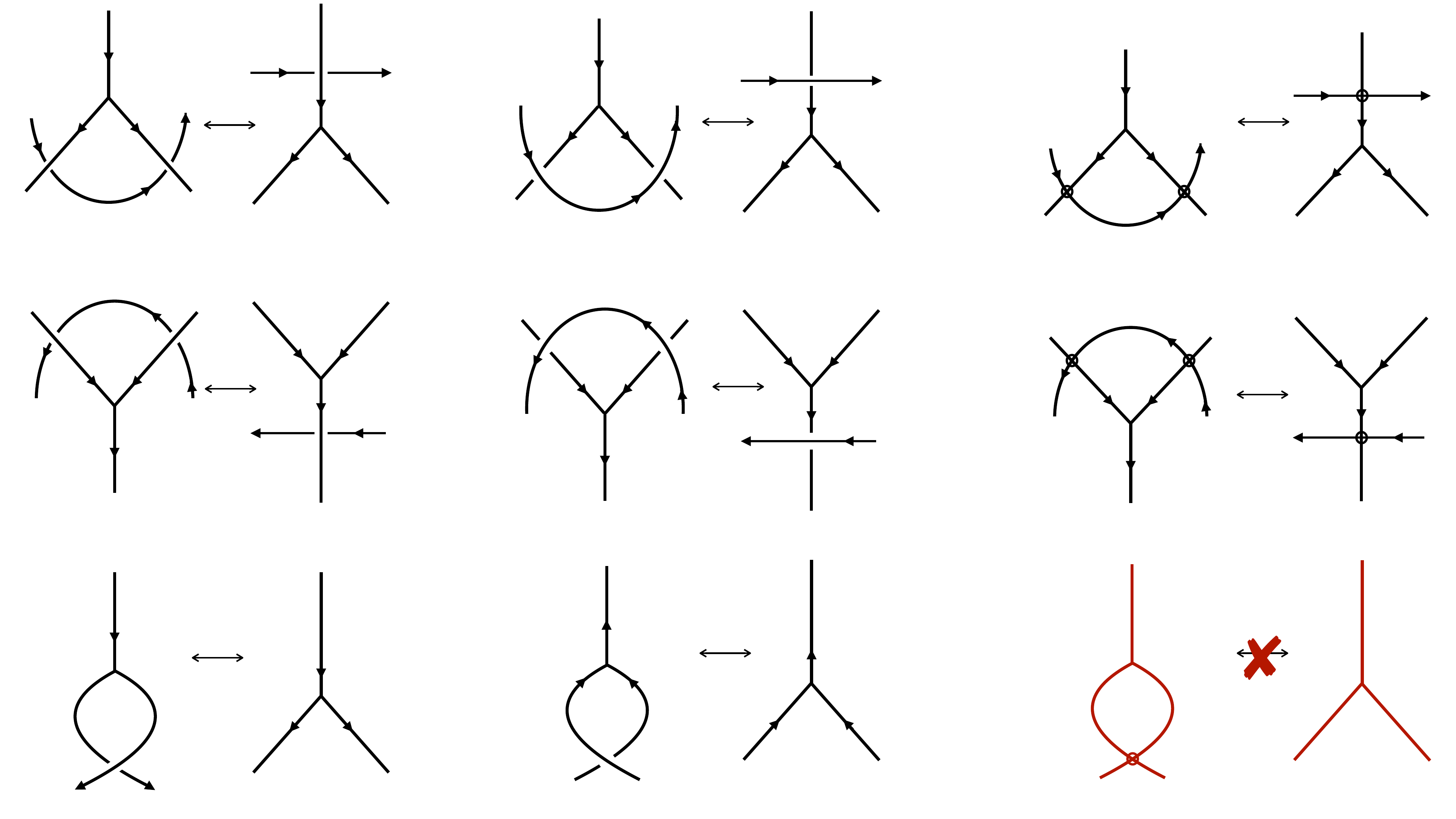}
    \put(-365,49){$R4$}
    \put(-217,49){$R4$}
    \put(-80,54){Forbidden}
    \put(-363,126){$R5$}
    \put(-215,126){$R5$}
    \put(-363,204){$R5$}
    \put(-217,204){$R5$}
     \put(-65,204){$vR5$}
    \put(-65,126){$vR5$}
    \caption{A subset of the virtual Y-oriented Reidemeister moves}
    \label{fig:vRMoves}
\end{figure}

 Graves, Nelson and the third author found a generating set of $Y$-oriented Reidemeister moves for $Y$-oriented trivalent spatial graphs \cite{GNT}. 
Flemming and Mellor studied virtual graphs with arbitrary degree and found a generating set of virtual graph Reidemeister moves in \cite{MT}. Combining the results of Graves-Nelson-T. and Flemming-Mellor, gives a complete list of Reidemeister moves for $Y$-oriented trivalent spatial graphs. 
Any set of generating Reidemeister moves for virtual knots (aka virtual Reidemeister moves) together with the moves in Figure \ref{fig:vRMoves} form a generating set of $Y$-oriented Reidemeister moves for virtual $Y$-oriented trivalent spatial graphs.

Notice, the virtual generalization of the $R4$ move is a forbidden move, and is not considered part of our virtual $Y$-oriented Reidemeister move list.

\section{From Tribrackets to Virtual \N Algebras}\label{sec:nalgs}

The new contributions of this paper are to synthesize two different generalizations of a \emph{tribracket}, namely the virtual tribracket and the \N Algebra, to create the virtual \N algebra. 
In this section, we explain each of these objects and their associated coloring invariants. We state the new definition of the virtual \N algebra along with its associated coloring invariant.

The term \emph{tribracket} was first coined by Nelson-Pico in \cite{NP} as an adaptation of the ternary quaisgroup defined in \cite{N} by Niebrzydowski.

\begin{defn}\label{defn:tri}(Niebrzydowski, Nelson-Pico)
    Let $X$ be a set. A \textbf{tribracket} on $X$ is a ternary operation $[-,-,-]:X\times X\times X\rightarrow X$ satisfying the following conditions.
    \begin{enumerate}
        \item  In the equation $[a,b,c]=d$, any three of the four elements in $\{a,b,c,d\}$ determine the fourth.
        \item For all $a,b,c,d\in X$ the following equalities are true.         \begin{align*}
             [a,b,[b,c,d]] &= [a,[a,b,c],[[a,b,c],c,d]]&\hfill(III.i) \\
             [[a,b,c],c,d] &= [[a,b,[b,c,d]], [b,c,d],d] &\hfill(III.ii)
         \end{align*}
    \end{enumerate}
\end{defn}

\begin{figure}[h]
    \centering
    \includegraphics[width=420\unitlength]{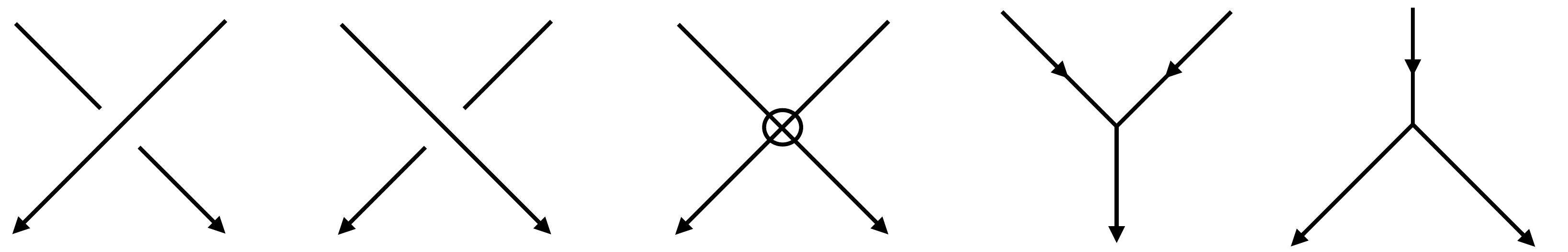}
    \put(-410,30){$a$}
    \put(-392,55){$b$}
    \put(-370,30){$c$}
    \put(-402,5){$[a,b,c]$}
    \put(-325,30){$a$}
    \put(-302,5){$b$}
    \put(-280,30){$c$}
    \put(-317,55){$[a,b,c]$}
    \put(-235,30){$a$}
    \put(-215,55){$b$}
    \put(-190,30){$c$}
    \put(-225,5){$\langle a,b,c\rangle $}
    \put(-140,30){$a$}
    \put(-125,55){$ab$}
    \put(-105,30){$b$}
    \put(-65,30){$a$}
    \put(-47,5){$ab$}
    \put(-25,30){$b$}
    \caption{A graphical interpretation for coloring instructions of planar complements of diagrams with tribrackets, virtual tribrackets, and \N  algebras.}
    \label{fig:coloring_instructions}
\end{figure}

The axioms of a  tribracket are derived from the  Reidemeister moves and can be interpreted as instructions for how to color the planar complement of a  link diagram, see Figure \ref{fig:coloring_instructions}. \N  proved that the number of colorings of the planar complement of an oriented  link diagram by a tribracket is an integer-valued invariant for classical links \cite{N, NP}.

There are two natural extensions of classical link theory;  virtual knot/link theory and trivalent spatial graphs. Nelson and Pico extend the notion of a tribracket to virtual knots in the form of a \emph{virtual tribracket}.

\begin{defn}\label{def:vtri}(Nelson-Pico)
A \textbf{virtual tribracket} on a set X is two ternary operations $[-,-,-], \langle -,-,-\rangle:X\times X\times X\rightarrow X$ which satisfy the following rules.
\begin{enumerate}

\item Both $[-,-,-]$ and $\langle -,-,-\rangle$ are tribrackets on $X$.

\item For all $a,b,c,d\in X$, the following equalities are true.

\begin{enumerate}

\item $\langle a,\langle a,b,c\rangle,c\rangle = b$\hfill(vII)
\item $[\langle a,b,c\rangle,c,d] =\langle [a,b,\langle b,c,d\rangle], \langle b,c,d\rangle, d\rangle$\hfill(m.i)
\item $[a,b,\langle b,c,d\rangle]=\langle a,\langle a,b,c\rangle, [\langle a,b,c\rangle, c,d]\rangle $\hfill(m.ii)
\end{enumerate}

\end{enumerate}
\end{defn}

The axioms of a virtual tribracket are derived from the virtual Reidemeister moves and can be interpreted as instructions for how to color the planar complement of a virtual link diagram, see Figure \ref{fig:coloring_instructions}. 
The square bracket $``[-,-,-]"$ is used for classical crossings, and the angle bracket $``\langle -,-,-\rangle"$ is used for virtual crossings. 
Nelson-Pico proved that the number of colorings of the planar complement of an oriented virtual link diagram by a virtual tribracket is an integer-valued invariant for virtual links \cite{NP}.

 In another direction, Graves-Nelson-T. introduced the \emph{\N algebras} as a tool to study $Y$-oriented trivalent spatial graphs \cite{GNT}. In a similar manner to the tribracket coloring invariant, one can also color the planar complement of a $Y$-oriented trivalent spatial graph with \N algebras to get an invariant of $Y$-oriented trivalent spatial graphs, see Figure \ref{fig:coloring_instructions}.

\begin{defn}\label{def:Nalg}(Graves-Nelson-T.)
Let $X$ be a set with a ternary operation $[-,-,-]:X\times X\times X\rightarrow X$ and a  product structure $a,b\mapsto ab$. $X$ is a \textbf{\N algebra} if the following conditions are satisfied.
\begin{enumerate}
    \item Any two of the three $\{a,b,c\}$ in $ab=c$ determines the third.
    \item $[-,-,-]$ is a tribracket on $X$.
    \item For all $a,b,c,d\in X$,the following equalities are true.\begin{enumerate}
    \item $[a,ab,b]=ab$ \hfill (R4)
        \item $a[a,b,c]=[a,b,bc]$\hfill (R5.1)
        \item $[a,b,c]=[[a,b,bc],bc,c]$ \hfill (R5.2)
        \item $[a,b,c]c=[ab,b,c]$\hfill (R5.3)
        \item $[a,b,c]=[a,ab,[ab,b,c]]$\hfill (R5.4)
        
    \end{enumerate}

\end{enumerate}
\end{defn}

The definition shown here is a slight modification of the definition given by Graves-Nelson-T. in \cite{GNT}.
 In the original definition of the \N algebra, the product could be a \emph{partially defined} binary operation. This leads to a welldefinedness issue for the coloring invariant, which we address in Section \ref{sec:welld}. The definition above requires the product to be fully defined on $X\times X$, and this leads to a well defined coloring invariant.

 The name \textit{algebra} is slightly misleading. The \N algebras are not algebras in the sense of a module with multiplication. But instead, these  \N algebras are \emph{like} an algebra of the  tribracket by adding a multiplication to the existing tribracket structure.

In the definition of a \N algebra, the relations $R4$ and $R5.1$, $R5.2$, $R5.3$, $R5.4$ are derived from the $R4$ and $R5$ $Y$-oriented Reidemeister moves. Coloring the planar complement of the arcs of the graph using the coloring rules in Figure \ref{fig:coloring_instructions} gives  the relations $R4$, $R5.1$, $R5.2$, $R5.3$, and $R5.4$. For more details, see Section 3 of \cite{GNT}.

In this paper, we combine the ideas of virtual tribrackets and \N algebras to get a \emph{virtual \N algebra}. We can color the regions of the planar complement of virtual $Y$-oriented trivalent spatial graph to obtain an invariant that generalizes both the tribracket and \N algebra color invariants.

\begin{defn}\label{def:vnalg}
A \textbf{virtual \N algebra} is a set $X$ with virtual tribracket given by  $[-,-,-]$ and $\langle -,-,-\rangle$, a  product $\cdot$ which maps $a,b\mapsto ab$, and satisfies the following conditions.
\begin{enumerate}
\item $[-,-,-]$ and $\langle -,-,-\rangle$ are a virtual tribracket as in Definition \ref{def:vtri}.

\item $X$ together with $[-,-,-]$ and the product $\cdot$ is a \N algebra as in Definition \ref{def:Nalg}.

\item For all $a,b,c\in X $,the following equalities are true. \begin{enumerate}
\item $a \langle a,b,c\rangle= \langle a,b,bc\rangle$ \hfill (vR5.1)
\item $\langle \langle a,b,bc\rangle, bc,c\rangle =\langle a,b,c\rangle$ \hfill (vR5.2)
\item $\langle a,b,c\rangle c=\langle ab,b,c\rangle$\hfill (vR5.3)
\item $\langle a,ab, \langle ab,b,c\rangle\rangle=\langle a,b,c\rangle$ \hfill (vR5.4)

\end{enumerate}

\end{enumerate}

\end{defn}

\noindent \textbf{Example 3.5.} The following brackets and product structure form a virtual \N algebras on the set $X=\{1,2,3\}$. The tribracket and multiplication structures are given in tables. For example, to determine $[a,b,c]$, take the the $b,c$-entry of the $a$'th matrix in $[-,-,-]$.

\begin{align*}[-,-,-]:&\left[ \begin{array}{rrr}1&2&3\\3&1&2\\2&3&1 \end{array}\right],\left[ \begin{array}{rrr}2&3&1\\1&2&3\\3&1&2 \end{array}\right],\left[ \begin{array}{rrr}3&1&2\\2&3&1\\1&2&3 \end{array}\right] &\cdot:\left[ \begin{array}{rrr}1&3&2\\3&2&1\\2&1&3 \end{array}\right]\\
\langle -,-,-\rangle:&\left[ \begin{array}{rrr}1&2&3\\3&1&2\\2&3&1 \end{array}\right],\left[ \begin{array}{rrr}2&3&1\\1&2&3\\3&1&2 \end{array}\right],\left[ \begin{array}{rrr}3&1&2\\2&3&1\\1&2&3 \end{array}\right]
\end{align*}

Like the virtual tribrackets in \cite{NP}, and \N algebras in \cite{GNT}, the equations in the definition of the virtual \N algebra are derived from coloring the planar complements of virtual $Y$-oriented Reidemseister moves. Figure \ref{fig:coloring_instructions} shows how to color regions of the planar complements of graphs using the tribracket and product structure of a virtual \N algebra.

As proved by Nelson-Pico, the number of colorings of a virtual knot by a virtual tribracket is an invariant of virtual knots \cite{NP}. That is, the axioms of virtual tribrackets are derived from, and hence preserve, the virtual Reidemeister moves.
Similarly, Graves-Nelson-T. proved that the number of \N algebra colorings of the planar complement of a $Y$-oriented trivalent spatial graph is an invariant of $Y$-oriented trivalent spatial graphs\cite{GNT}.
That is, all of the (classical) $Y$-oriented Reidemeister moves and the classical Reidemeister moves are preserved by the axioms of a \N algebra. 
Virtual \N algebras are \N algebras with an extra virtual tribracket.
Combining the two results for tribrackets and \N algebras, to show that the number of colorings of a virtual $Y$-oriented trivalent spatial graph by a virtual \N algebra is an invariant, we only need to check that the two virtual $Y$-oriented Reidemeister moves ($vR5$) are preserved under these colorings. 

\begin{figure}
    \centering
    \includegraphics[width=300\unitlength]{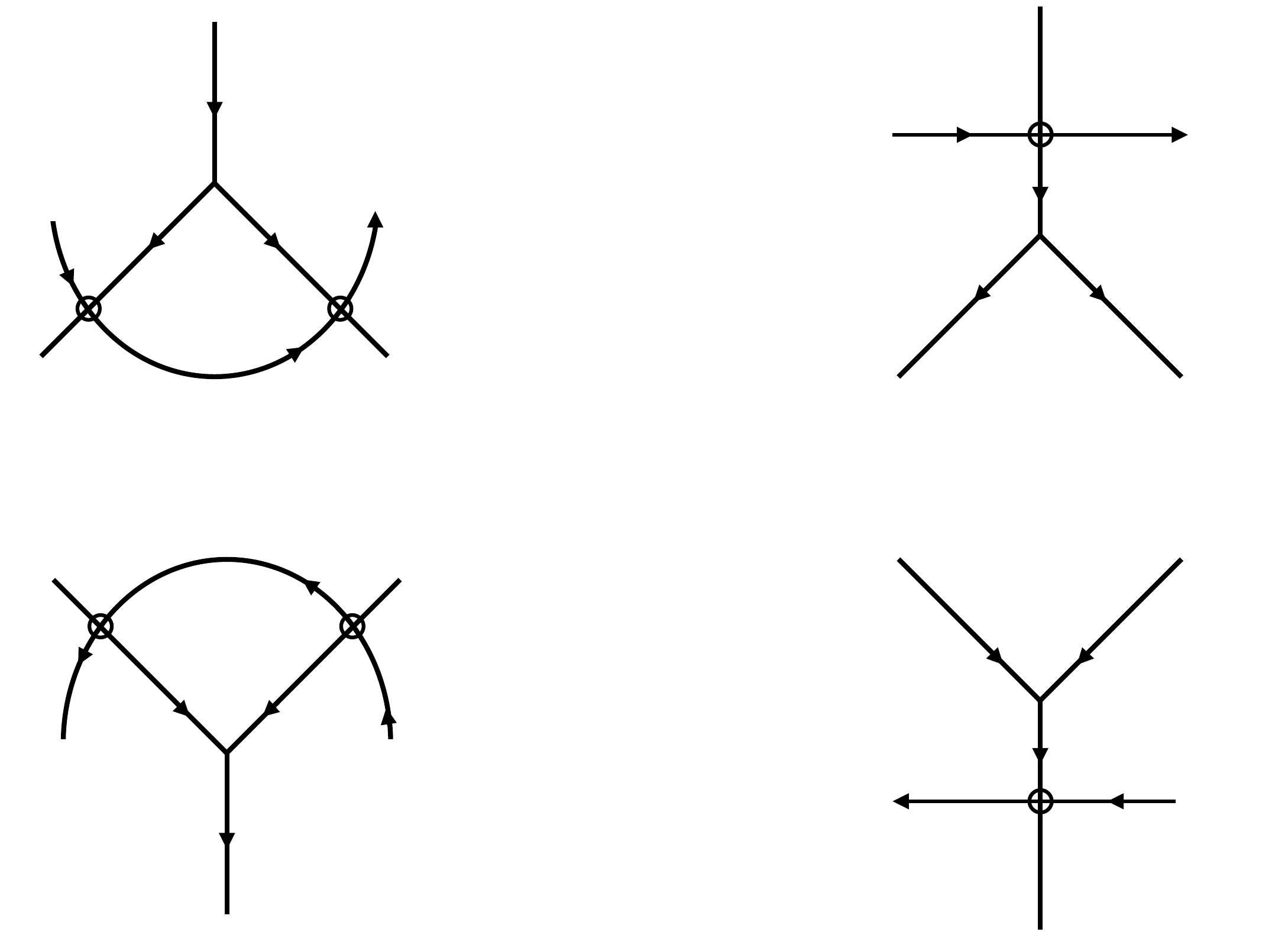}
    \put(-300,150){$a$}
    \put(-275,175){$b$}
    \put(-230,175){$c$}
    \put(-255,145){$bc$}
    \put(-268,123){\textcolor{blue}{$\langle a,b,bc\rangle$}}
    \put(-215,150){\textcolor{red}{$\langle \langle a,b,bc\rangle,bc,c\rangle$}}
    \put(-170,210){$vR5.1$}
    \put(-170,200){$\longleftrightarrow$}
    \put(-170,190){$vR5.2$}
    \put(-85,165){$a$}
    \put(-85,205){$b$}
    \put(-40,205){$c$}
    \put(-40,165){\textcolor{red}{$\langle a,b,c\rangle$}}
    \put(-75,125){\textcolor{blue}{$a\langle a,b,c\rangle$}}
    \put(-300,70){$a$}
    \put(-255,100){$ab$}
    \put(-202,70){$b$}
    \put(-265,70){\textcolor{blue}{$\langle ab,b,c\rangle$}}
    \put(-355,37){\textcolor{blue}{$\langle a,ab,\langle ab,b,c\rangle \rangle c$}}
    \put(-272,38){\begin{tikzpicture}
    \draw[->,blue] (0,0) to[out=0,in=205] (.8,.7);
    \end{tikzpicture}}
    \put(-225,20){$c$}
    \put(-322,20){\textcolor{red}{$\langle a,ab,\langle ab,b,c\rangle \rangle$}}
    \put(-170,40){$vR5.3$}
    \put(-170,30){$\longleftrightarrow$}
    \put(-170,20){$vR5.4$}
    \put(-95,15){\textcolor{red}{$\langle a,b,c\rangle $}}
    \put(-40,15){$c$}
    \put(-40,45){$b$}
    \put(-80,45){$a$}
     \put(-63,80){$ab$}
    \caption{Two virtual $Y$-oriented Reidemeister moves which describe the virtual \N algebra relations $vR5.1$, $vR5.2$, $vR5.3$, and $vR5.4$.}
    \label{fig:vGraphMoves}
\end{figure}

Figure \ref{fig:vGraphMoves} shows the two orientations of the virtual $Y$-oriented Reidemeister move $vR5$ colored according to the rules of the product $\cdot$ and tribracket $\langle -,-,-\rangle$. 
In the definition of a virtual \N algebra, axioms $vR5.1$ and $vR5.2$ can be seen by equating the two blue and two red quantities in the top diagrams in Figure \ref{fig:vGraphMoves}. Axiom $vR5.4$ can be seen by equating the two red quantities in the bottom diagrams in Figure \ref{fig:vGraphMoves}. By equating the blue quantities in the bottom diagram of Figure \ref{fig:vGraphMoves}, we get the equation $\langle a,ab,\langle ab,b,c\rangle\rangle c= \langle ab,b,c\rangle$. By substituting $\langle a,ab,\langle ab,b,c\rangle\rangle= \langle a,b,c\rangle$ from axiom $vR5.4$, we get $vR5.3$.

By construction, the following theorem is true.

\begin{thm}
Let $X$ be a virtual \N algebra and $\Gamma$ a virtual $Y$-oriented trivalent spatial graph. The number $\Phi_X(\Gamma)$ of $X$-colorings of the planar complement of $\Gamma$ is invariant under the virtual $Y$-oriented Reidemeister moves, and is therefore an integer-valued invariant of virtual $Y$-oriented trivalent spatial graphs.

\end{thm}

\noindent{\textbf{Example 3.7.}
    Using $X=\{1,2,3\}$ and the virtual \N algebra from Example 3.5, Figure \ref{fig:2difcolorings} shows two virtual Y-oriented trivalent spatial graphs that are distinguished by the invariant $\Phi_X$.  For the graph $K_1$, choosing the labels $a,b$, and $c$ determines the colorings for every other region and is limited by one equation $a\langle a, ac,a\rangle=\langle a,ac,c\rangle c.$ There are three possible choices for $a,b$, and $c$ from $X$, which makes $\Phi_X(K_1)=3$.

    For the graph $K_2$, choosing the labels $a,b,c,d$ and $e$ determines the coloring of the entire planar complement and is limited by the following seven equations. 
    
    \begin{align*}
        bc&=b((bc)[bc,c,a])&\langle b,(bc)[bc,c,a],[bc,c,a]a\rangle&=c\\
        [bc,c,a]&=((bc([bc,c,a])([bc,c,a]a)& [[bc,c,a]a,a\langle a,d,e \rangle,e]&=d\\
        a&=([bc,c,a]a)(a\langle a,d,e\rangle) & [c,a,d]&=[bc,c,a]a\\
        \langle a,d,e\rangle&=(a\langle a,d,e\rangle)e \\     
    \end{align*}
There are 9 possible choices for the labels $a,b,c,d$ and $e$ from $X$, which makes $\Phi_X(K_2)=9$. Since $\Phi_X(K_1)\neq \Phi_X(K_2)$, then $K_1\neq K_2$.

    \begin{figure}
        \centering
        \includegraphics[scale=.4]{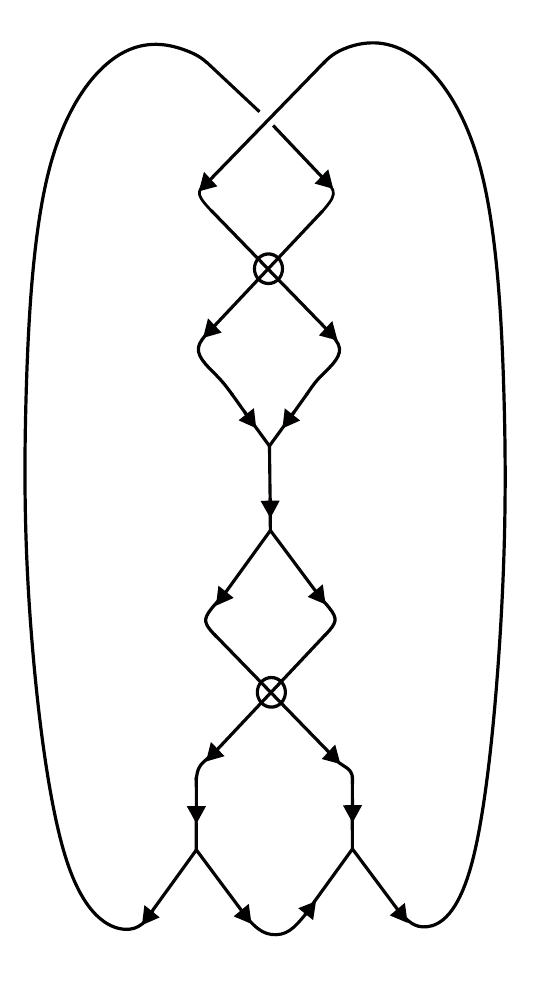}
        \includegraphics[scale=.4]{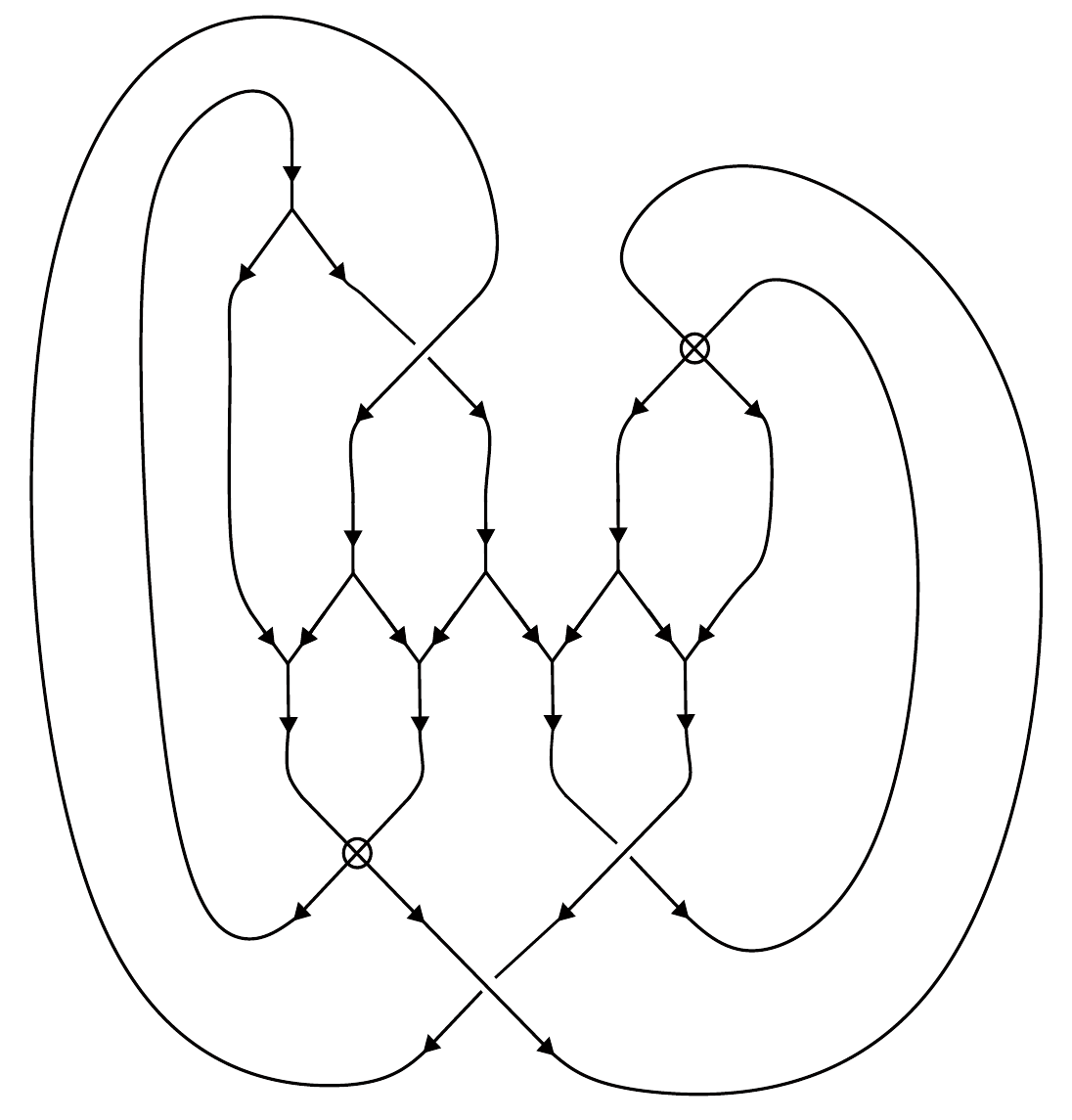}
        \put(-320,0){$K_1$}
        \put(-10,0){$K_2$}
        \put(-295,165){$a$}
        \put(-273,175){$b$}
        \put(-255,165){$c$}
        \put(-105,150){$a$}
        \put(-165,185){$b$}
        \put(-135,170){$c$}
        \put(-75,165){$d$}
        \put(-63,150){$e$}
        \caption{Two virtual Y-oriented trivalent spatial graphs $K_1$ and $K_2$. Three planar regions in the complement of $K_1$ have been labeled with $a,b$, and $c$. Five planar regions in the complement of $K_2$ have been labeled  by $a,b,c,d$, and $e$.}
        \label{fig:2difcolorings}
    \end{figure}

\section{Partial Products and Welldefinedness}\label{sec:welld}

In their original paper, Graves-Nelson-T. defined a \N algebra  as in Definition \ref{def:Nalg} but where the product structure need only be partially defined.  That is, if the product is written as a table, some entries of the table are left blank, as shown below. We will call a \N algebra with a partially defined product structure a \emph{partially defined \N algebra}. 
\[\left[ \begin{array}{rrr}1&-&2\\3&2&-\\-&-&3 \end{array}\right]\]

Using a partially defined \N algebra to color the planar complement of a graph restricts the colorings to only include those that give rise to defined products. For example, in Figure \ref{fig:undefinedtrivlaentcoloring}, if $ab$ is not defined, then $a$ and $b$ cannot color regions on either side of a trivalent edge.
 
\begin{figure}
\begin{tikzpicture} [scale=.75]
\draw[thick](-1,1)--(0,0);
\draw[thick,->](1,1)--(0,0)--(0,-1);
\node  at(-.75,0) {$a$};
\node  at(.75,0) {$b$};
\node at (0,.75){$ab$};
\node  at(0,1.15) {{\Tiny \textit{undefined}}};
\begin{scope}[xshift=4cm]
\draw[thick,->](0,1)--(0,0);
\draw[thick](-1,-1)--(0,0)--(1,-1);
\node  at(-.75,.1) {$a$};
\node  at(.75,.1) {$b$};
\node at (0,-.75){$ab$};
\node  at(0,-1.15) {{\Tiny \textit{undefined}}};
\end{scope}
\end{tikzpicture}
\caption{Partial products restrict possible colorings around a trivalent vertex.}
\label{fig:undefinedtrivlaentcoloring}
\end{figure}

There are several benefits of using a partially defined \N algebra. The partial product restricts the number of colorings of a graph, making the associated coloring invariant easier to compute. 
Also, there are not many examples of fully defined \N algebras, as will be described in Section \ref{sec:comp}. Allowing a partial product gives a much larger family of \N algebra coloring invariants. 
However, the use of a partially defined product structure can lead to a welldefinedness issue.

Suppose $X$ is a partially defined \N algebra where $b,c\in X$ but $bc$ is undefined. Consider the colorings in Figure \ref{fig:GraphMoves} before and after a Reidemeister move. The coloring on the right does not use $bc$, the undefined product, where the coloring on the left does use the product $bc$.
Therefore the coloring on the right is a valid coloring and the coloring on the left is not. 
This shows that coloring the planar complement by a partially defined \N algebra does not respect the Reidemeister moves and is not an invariant of the graph.

\begin{figure}
    \centering
    \includegraphics[width=300\unitlength]{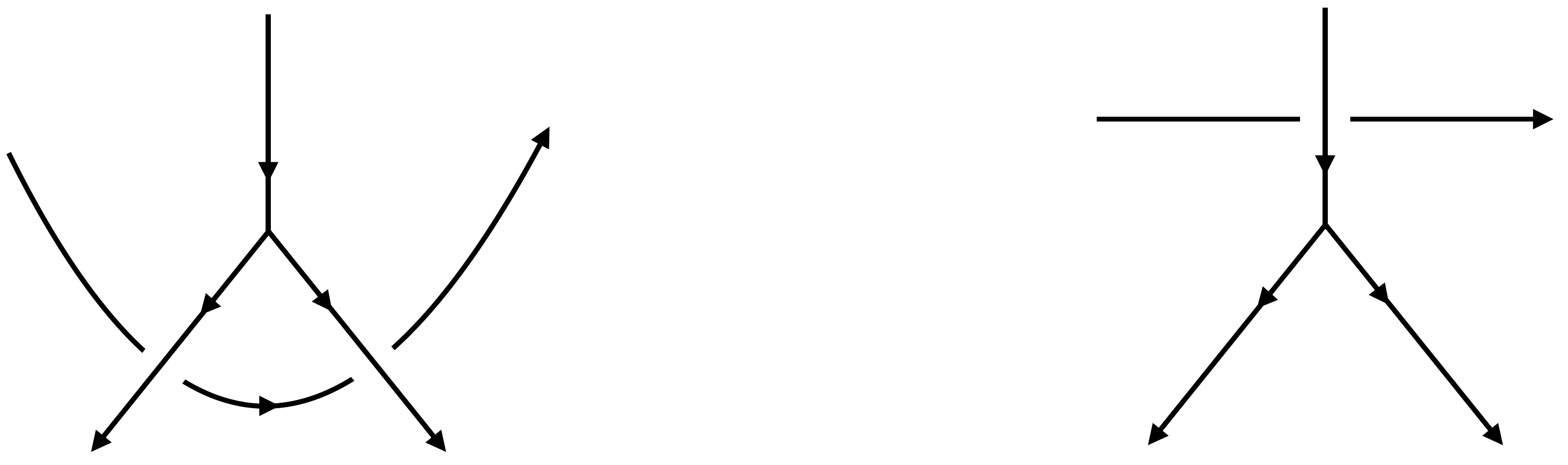}
    \put(-300,20){$a$}
    \put(-275,45){$b$}
    \put(-230,45){$c$}
    \put(-255,20){\textcolor{red}{$bc$}}
    \put(-268,-5){\textcolor{red}{$[ a,b,bc]$}}
    \put(-215,20){\textcolor{red}{$[ [a,b,bc],bc,c]$}}
    \put(-150,55){$\longleftrightarrow$}
    \put(-75,35){$a$}
    \put(-75,75){$b$}
    \put(-25,75){$c$}
    \put(-30,35){$[ a,b,c]$}
    \put(-65,-5){$a[ a,b,c]$}
    \caption{A $Y$-oriented Reidemeister can introduce an undefined product.}
    \label{fig:GraphMoves}
\end{figure}

We give a solution to this welldefinedness issue.
Keeping the definition of the \N algebra to include the partial products, we can instead define the coloring invariant associated to the partially defined  \N algebra to count only \emph{admissible colorings} of a graph. An admissible coloring is a coloring of the planar complement where an undefined product is never induced by a Reidmemeister move. This would of course lead to a well defined graph invariant, but it is very difficult to compute as it is difficult to prove a coloring is admissible, and would be very graph dependent.

To make this coloring invariant computable, we can add an axiom to the definition of an \N algebra that will force every coloring by that \N algebra to be admissible.
At first glance, it is tempting to add to the definition of \N algebra that for every axiom, the left side of the equation is defined if and only the right side is defined. 
For example, $(3)(a)$ of Definition \ref{def:Nalg} states that \[ [a,ab,b]=ab.\]

In this example, it is true that $[a,ab,b]$ is defined if and only if $ab$ is defined. However, in axiom $(3)(c)$, \[ [a,b,c]=[[a,b,bc],bc,c],\]
 $[a,b,c]$ is always defined for all $a,b,c$, so $bc$ must always be defined. This would imply that the product structure is fully defined.

Instead, we focus on the the Reidemeister moves to derive the axiom. The axiom needs to ensure that the seemingly admissible coloring on the right side of Figure \ref{fig:GraphMoves} never occurs. So if $bc$ is undefined, then $a$ and $a[a,b,c]$ cannot color regions on either side of a trivalent vertex. One way to ensure this is if $a[a,b,c]$, the product of $a$ and $[a,b,c]$, is also undefined.
Taking into account the same Reidemeister move with the negative crossings, the new axiom is:
\[\text{For all $b$ and $c$, } bc \text{ is defined if and only if }a[a,b,c] \text{ and } [b,c,a]a \text{ are defined for all $a$}.\]

The corrected definition of a partially defined \N algebra with the new axiom is stated below.  

\begin{defn}\label{def:wdNalg}
Let $X$ be a set with a ternary operation $[-,-,-]:X\times X\times X\rightarrow X$ and a partial product $a,b\mapsto ab$ structure. $X$ is a \textbf{partially defined \N algebra} if the following conditions are satisfied.
\begin{enumerate}
    \item If the product $ab$ is defined, any two of the three $\{a,b,c\}$ in $ab=c$ determines the third.
    \item $[-,-,-]$ is a tribracket on $X$.
    \item For all $a,b,c,d\in X$, if both sides of the equations are defined, then the following equalities are true.\begin{enumerate}
    \item $[a,ab,b]=ab$ \hfill (R4)
        \item $a[a,b,c]=[a,b,bc]$\hfill (R5.1)
        \item $[a,b,c]=[[a,b,bc],bc,c]$ \hfill (R5.2)
        \item $[a,b,c]c=[ab,b,c]$\hfill (R5.3)
        \item $[a,b,c]=[a,ab,[ab,b,c]]$\hfill (R5.4)
        
    \end{enumerate}
    
\item  For all $b$ and $c$, $bc$ is defined if and only if $a[a,b,c]$ and $[b,c,a]a$ are defined for all $a$.
\end{enumerate}
\end{defn}

The addition of the new axiom, axiom (4), ensures that any coloring by a partially defined \N algebra will be an admissible coloring and respect the Reidemeister moves. 

\begin{cor}
Any coloring of the planar complement of a Y-oriented spatial graph by a partially defined \N algebra is an admissible coloring. The total number of admissible colorings by a partially defined \N algebra is a graph invariant.
\end{cor}

Applying the same reasoning to the virtual analogue of the Reidemeister move shown in Figure \ref{sec:welld}, we can add the same axiom for the $\langle -,-,-\rangle$ bracket. This will ensure that every coloring by a partially defined virtual \N algebra is an admissible coloring.

\begin{defn}\label{def:wdvnalg}
A \textbf{partially defined virtual \N algebra} is a set $X$ with virtual tribracket given by  $[-,-,-]$ and $\langle -,-,-\rangle$, a partial product $\cdot$ which maps $a,b\mapsto ab$, and satisfies the following conditions.
\begin{enumerate}
\item $[-,-,-]$ and $\langle -,-,-\rangle$ are a virtual tribracket as in Definition \ref{def:vtri}.

\item $X$ together with $[-,-,-]$ and the product $\cdot$ is a partially defined \N algebra as in Definition \ref{def:wdNalg}.

\item For all $a,b,c\in X $, if both sides of the equations are defined, then the following equalities are true, \begin{enumerate}
\item $a \langle a,b,c\rangle= \langle a,b,bc\rangle$ \hfill (vR5.1)
\item $\langle \langle a,b,bc\rangle, bc,c\rangle =\langle a,b,c\rangle$ \hfill (vR5.2)
\item $\langle a,b,c\rangle c=\langle ab,b,c\rangle$\hfill (vR5.3)
\item $\langle a,ab, \langle ab,b,c\rangle\rangle=\langle a,b,c\rangle$ \hfill (vR5.4)

\end{enumerate}
\item For all $b$ and $c$, $bc$ is defined if and only if $a\langle a,b,c\rangle$ and $\langle b,c,a\rangle a$ are defined for all $a$.
\end{enumerate}

\end{defn}

\begin{cor}
Any coloring of the planar complement of a virtual Y-oriented spatial graph by a partially defined virtual \N algebra is an admissible coloring. The total number of admissible colorings by a partially defined virtual \N algebra is a graph invariant.
\end{cor}

\begin{ex} The set $X=\{1,2,3,4\}$ with the following tribracket and multiplicative structures forms a partially defined virtual \N algebra.

\begin{align*}[-,-,-]:&\left[ \begin{array}{rrrr} 1& 2& 3& 4\\4& 1& 2& 3\\3& 4& 1& 2\\ 2& 3& 4& 1\end{array}\right],
\left[ \begin{array}{rrrr}2& 3& 4& 1\\1 &2 &3 &4\\4& 1& 2& 3\\ 3 &4& 1& 2 \end{array}\right],
\left[ \begin{array}{rrrr}3& 4& 1& 2\\2& 3& 4& 1\\1 &2 &3& 4\\4& 1& 2& 3 \end{array}\right],
\left[
\begin{array}{rrrr}4& 1& 2& 3\\3& 4& 1& 2\\2& 3& 4& 1\\ 1& 2& 3& 4 \end{array}\right]\\ 
\langle -,-,-\rangle :&\left[ \begin{array}{rrrr}1 &2 &3 &4\\4& 1& 2& 3\\3& 4& 1& 2\\2& 3& 4& 1 \end{array}\right],
\left[ \begin{array}{rrrr}2& 3& 4& 1\\1& 2& 3& 4\\4& 1& 2& 3\\3& 4& 1& 2 \end{array}\right],
\left[ \begin{array}{rrrr}3& 4& 1& 2\\2& 3& 4& 1\\1& 2& 3& 4\\4& 1& 2& 3 \end{array}\right],
\left[
\begin{array}{rrrr}4 &1& 2& 3\\3& 4& 1& 2\\2& 3& 4& 1\\1& 2& 3& 4 \end{array}\right]\\
\cdot:&\left[\begin{array}{rrrr}1& -& -& -\\-& 2& -& -\\-& -& 3& -\\-& -& -& 4 \end{array}\right]
\end{align*}

\end{ex}

\section{Computational Results}\label{sec:comp}

How many different coloring invariants are there? That is, how many tribrackets, \N algebras, and virtual \N algebras  are there?
Using Latin squares and Latin cubes, we exhaustively generated lists of all tribrackets, \N algebras and virtual \N algebras on a set with 3 elements, and a set with 4 elements.
All of the code to verify our results and to generate usable lists of tribrackets and \N algebras can be found on GitHub in \cite{GIT}.

\subsection{Latin Squares and Cubes.}
A \emph{Latin square of order $n$} is an $n\times n$ array filled with $n$ different symbols, each occurring exactly once in each row and exactly once in each column. 
Let $X$ be a set with $n$ elements and $\mathcal{N}$ a \N algebra on $X$. 
Condition (1) in Definition \ref{def:Nalg} of a \N algebra sates for any two of the three $\{a,b,c\}\subseteq X$ in the equation $a\cdot b=c$ determines the third. This condition requires that any complete multiplication table for the product operation must be a Latin square. 
On the other hand, for a partially defined product, any two of the three $\{a,b,c\}\subseteq X$ in the equation $a\cdot b=c$ determines the third only implies no row or column of the multiplication table can contain a repeated entry. 
From this definition, a partially defined product does not need to be completable to a Latin square.

A \emph{Latin cube of order $n$} is an $n\times n\times n$ cube array filled with $n$ different symbols, each occurring exactly once in each row, exactly once in each column, and exactly once in each third dimension. Condition (1) in Definition \ref{defn:tri} of a tribracket requires that any multiplication table for the tribracket operation must be a Latin cube of order $n$.

The number of Latin squares of order $n$ has a closed formula, though not easily computable, and has been been computed up to $n=10$ in the \emph{Online Encyclopedia of Integer Sequences} (OEIS) \cite{OEIS} as sequence A002860. 
There are 12 Latin squares order three, 576 Latin squares order four,  and 161,280 Latin squares order five.
There are 24  Latin cubes or order three, 55,296  Latin cubes or order four, and 2,781,803,520  Latin cubes of order five, as verified in sequence A098679 in the OEIS.  

\subsection{Tribrackets}\label{sec:tribracket}
To generate tribrackets, we first generated all Latin cubes of order 3 and 4, and then exhaustively checked which cubes satisfied the tribracket axioms. 
 Our computations have revealed that there exist 12 distinct tribrackets order 3 and 168 distinct tribrackets order 4. We guess that there are 4680 tribrackets of order 5, and in general $n(n^{(n-1)}-1)\frac{(n-2)!}{(n-1)}$ tribrackets of order $n$, based on the sequence A059522 from OEIS.  
 Our computational approach was not efficient enough to exhaustively generate and search through the 2.7 billion Latin Cubes of order 5, but we encourage the reader to be more clever than we were.

By considering all pairs of tribrackets and checking to see which pairings satisfy the virtual tribracket axioms, we exhaustively computed all virtual tribrackets of order 3 and 4.
We found all 24 distinct virtual Tribrackets order 3, and 1080 distinct virtual Tribrackets of order 4. 

\subsection{\N Algebras}\label{sec:niebrzydowski}
Compared to the plethora of tribrackets, \N algebras are quite scarce.
By exhaustively searching all pairings of a tribracket with a Latin square, 
we demonstrated there is a unique \N algebra of order 3 (the one shown in Example 3.5), and no \N algebras of order 4.
We suspect this inductively implies there are no \N algebras of order $n>4$, though we cannot prove it.

Loosening the restrictions to allow for a partially defined product, we determined there are 6 partially defined \N algebras of order 3.
All instances used the same tribracket as in Example 3.5, and the 6 different partial products are listed below.
 It turns out that for order 3, every partial product square can be completed to a fully defined Latin square.

\[
\left[\begin{array}{rrr}-& -& 2\\3& -& -\\-& 1& -
\end{array}\right]
,\left[\begin{array}{rrr}1& 3& -\\-& 2& 1\\ 2& -& 3 \end{array}\right]
,\left[\begin{array}{rrr}-& 3& -\\ -& -& 1\\2& -& - \end{array}\right],\] \[
\left[\begin{array}{rrr}1& -& 2\\ 3& 2& -\\-& 1& 3 \end{array}\right]
,\left[\begin{array}{rrr}1& -& -\\ -& 2& -\\-& -& 3 \end{array}\right]
,\left[\begin{array}{rrr}-& 3& 2\\ 3& -& 1\\2& 1& - \end{array}\right]
\]

We also found 28 partially defined \N algebras of order 4. We only searched for partial squares that can be completed to a fully defined Latin square. It is possible that there are examples that we did not find that have a product which does not complete to a Latin square.

\subsection{Virtual \N Algebras}\label{sec:virtual-niebrzydowski}

Similarly to the \N algebras, we found there to be a unique virtual \N algebra of order 3 (as shown in Example 3.5), and since there is no \N algebra of order 4, there is no virtual \N algebra of order 4.

There are exactly 6 partially defined virtual \N algebras of order 3, all which come from the 6 examples of the partially defined \N algebras listed above.

Lastly, we found 52 partially defined virtual \N algebras of order 4. We only searched for partial products that can be completed to a fully defined Latin cube, so there could potentially be more examples than found here.
The order $4$ partial squares have an entry for every row and for every column, and fit into one of the patterns shown below.
\[
\left[\begin{array}{rrrr}*& -& -& -\\-& *& -& -\\-& -& *& -\\-& -& -& * \end{array}\right]
,\left[\begin{array}{rrrr}-& -& *& -\\-& -& -& *\\ *& -& -& -\\-& *& -& - \end{array}\right]
,\left[\begin{array}{rrrr}-& *& -& -\\ *& -& -& -\\-& -& -& *\\-& -& *& - \end{array}\right]
,\left[\begin{array}{rrrr}-& -& -& *\\-& -& *& -\\-& *& -& -\\ *& -& -& - \end{array}\right]
\]

\subsection{Small mysteries for the interested reader.}
Our computational findings seem to beg more questions than they answer. We encourage the reader to consider the following questions. Why are there no \N algebras of order 4? 
Where is the contradiction in the axioms for $n=4$ that is not a contradiction for $n=3$? Does this contradiction apply to all $n>4$? Can one inductively prove that if there are no \N algebras of order $4$, then there are no \N algebras of order $n>4$? For the virtual \N algebras of order 4, why are the missing entries in the product squares symmetric?

\bibliography{virtualNalgebras}{}
\bibliographystyle{plain}

\end{document}